# AVERAGE INTERPOLATING WAVELETS ON POINT CLOUDS AND GRAPHS

## RAIF M. RUSTAMOV

ABSTRACT. We introduce a new wavelet transform suitable for analyzing functions on point clouds and graphs. Our construction is based on a generalization of the average interpolating refinement scheme of Donoho. The most important ingredient of the original scheme that needs to be altered is the choice of the interpolant. Here, we define the interpolant as the minimizer of a smoothness functional, namely a generalization of the Laplacian energy, subject to the averaging constraints. In the continuous setting, we derive a formula for the optimal solution in terms of the poly-harmonic Green's function. The form of this solution is used to motivate our construction in the setting of graphs and point clouds. We highlight the empirical convergence of our refinement scheme and the potential applications of the resulting wavelet transform through experiments on a number of data stets.

## 1. INTRODUCTION

In many applications, the data naturally comes in the form of functions defined on non-standard domains such as point clouds and graphs. Of great practical interest are tools suitable for the processing of such functions, including denoising, compression, and learning. These tasks have been extensively studied in the setting of standard Euclidean domains, where multiscale methods have proven to be quite effective. One of the reasons for such effectiveness is, perhaps, the sheer quantity of various multiscale constructions available on standard domains. In contrast, on non-standard domains, there is a very limited choice of multiscale constructions; these include the Haar-like wavelets [17, 8], diffusion wavelets [5], lifting wavelets [10], and spectral wavelets [9].

The main goal of this paper is to increase the assortment of available choices by generalizing the average interpolating wavelet scheme of Donoho [7] to point clouds and graphs. Our interest in this particular scheme is due to the following two reasons. First, as explained in [14], when analyzing a function consisting purely of white noise, the resulting average interpolating wavelet coefficients are basically noise, but have roughly the same size at all scales and locations, and are roughly independent. These properties make the average interpolating wavelets especially useful for denoising. Second, the approach to regression via Haar wavelets proposed in [8], can be easily generalized to the average interpolating wavelets, making them suitable for machine learning problems.

The average interpolating scheme is based on predicting function averages on smaller children regions based on the given averages on the larger parent region and its neighbors. This is achieved by fitting an interpolant to the provided averages, and inferring the averages on the children regions from this interpolant. The

---

construction of the interpolant in the setting of the real line is straightforward – polynomials provide a natural choice. For example, if one bases the interpolation on the immediate neighbors of the parent region (dyadic interval) on the real line, there are three averages (two for neighbors, and one for the interval itself) to which one can fit a unique quadratic polynomial. If one wants to run the scheme with second order neighborhoods (neighbors and neighbors of neighbors), then there are five averages to which a fourth degree polynomial can be fit; similarly for higher order neighborhoods. Note that a polynomial of the same fixed degree is used throughout the construction.

These conveniences disappear in the case of manifolds, point clouds, and graphs. First, there is no straightforward way of defining polynomials on manifolds, which means that we need a new interpolation basis. Second, the number of neighborhoods at every region can be different, which at first glance means that an interpolation basis of different dimension will be needed at each region. However, if we were to use such order variable bases, the issue of stability would arise. Namely, when trying to predict the averages for the children of two neighboring regions, we may end up fitting very different interpolants and obtaining very disparate predicted averages at the children. This would make obtaining continuous wavelets impossible.

Thus, the main challenge in obtaining a useful wavelet construction on general domains is to design a stable interpolation scheme that can handle varying numbers of neighbor regions. An additional requirement is that if all of the averages at the parent region and its neighbors are the same, then the interpolant must be constant. This is needed in order to ensure that the wavelets have a vanishing moment. To overcome this challenge, we define the interpolant as the minimizer of a smoothness functional subject to the averaging constraints. In the continuous setting, we derive a formula for the optimal solution in terms of the poly-harmonic Green's function. The form of this solution is then used to motivate our construction in the setting of graphs and point clouds. The resulting discrete interpolant has an extremely simple form: it is a linear combination of low-frequency (i.e. corresponding to smaller eigenvalues) eigenfunctions of the graph Laplacian. Since there can be many linear combinations that satisfy the averaging constraints, we single out a linear combination among them by picking the one that minimizes a measure of smoothness, such as the spline/Laplacian energy.

This paper is organized as follows. Section 2 introduces an outline of our construction in the setting of compact Riemannian manifolds without boundary, the main goal being to provide an easy to follow presentation of the main ideas, and to motivate some of the choices made later in the discrete setting. Our presentation in this section is not meant to be rigorous (e.g. we do not know/prove that the refinement scheme converges), as we are truly interested in the ramifications in the discrete case. Section 3 provides the details and algorithms for the setting of weighted graphs and point clouds. Section 4 highlights the empirical convergence of our refinement scheme and the potential applications of the resulting wavelet transform through experiments on a number of data stets.

## 2. Motivation: Closed Manifolds

2.1. **Average Interpolating Wavelets.** In this subsection we review the main steps of the average interpolating wavelet construction of Donoho [7] with modifications as needed in the setting of manifolds; our exposition closely follows the outline of [14].

For a closed Riemannian manifold $M$, we first construct a hierarchical dyadic partition of $M$ into contiguous regions $\{R_{l,k}, l = 0, 1, ..., \text{and } k = 1, ..., 2^l\}$. We start by setting $R_{0,1} = M$, and then partition $R_{0,1}$ using an appropriate procedure into two regions $R_{1,1}$ and $R_{1,2}$. This process is repeated on each of the obtained regions in recursion, by partitioning each region $R_{l,k}$ into two children regions $R_{l+1,2k-1}$ and $R_{l+1,2k}$. As a result, we obtain infinitely many regions organized into a binary tree. We assume that the sizes (volumes, diameters) of regions at the same level $l$ are comparable and converge to zero with an increasing $l$. For a given region $R_{l,k}$ at level $l$, we use $Vol(R_{l,k})$ to denote its volume. At each level $l$, we can compute the average volume of all of the regions at that level; we denote this average by $VolAv_l$. For region $R_{l,k}$ at level $l$, we can consider its neighboring regions, i.e. regions at the same level $l$ that have a common boundary with $R_{l,k}$; we denote these regions by $R_{l,(k,r)}$ with $r = 1, ..., m$, where $m$ is the number of neighbors, including the region itself; i.e. for convenience, we set $(k, 1) = k$.

Given an integrable real-valued function $f : M \to \mathbb{R}$, consider a pyramid of averages

$$\beta_{l,k} = Ave\{f|R_{l,k}\}$$

of function $f$ over the regions $R_{l,k}$. Since we have $R_{l,k} = R_{l+1,2k-1} \cup R_{l+1,2k}$, it follows that the pyramid is redundant with the following two scale relation

$$(2.1) \qquad \beta_{l,k} = \frac{\beta_{l+1,2k-1} Vol(R_{l+1,2k-1}) + \beta_{l+1,2k} Vol(R_{l+1,2k})}{Vol(R_{l+1,2k-1}) + Vol(R_{l+1,2k})}.$$

In the average interpolating wavelet construction, one tries to predict the averages of $f$ at the level $l + 1$ from the averages at the level $l$, and the discrepancies between predicted and actual averages are recorded as the wavelet coefficients.

To be more precise, assume that we are given the averages $\beta_{l,k}$ where $k = 1, ..., 2^l$ and we would like to produce the predicted averages at the next level, $\tilde{\beta}_{l+1,k}$, $k = 1, ..., 2^{l+1}$. For a fixed region $R_{l,k}$, consider all of its neighboring regions $R_{l,(k,r)}$. Fit a more or less regular function $\pi_{l,k} : M \to \mathbb{R}$ to the average values provided over these regions $R_{l,(k,r)}$; namely, construct a function $\pi_{l,k}$ which satisfies

$$(2.2) \qquad Ave\{\pi_{l,k}|R_{l,(k,r)}\} = \beta_{l,(k,r)}, \, r = 1, ..., m,$$

where $m$ is the number of neighbors of our region (remember, the region itself is also included in this count, and $R_{l,(k,1)} = R_{l,k}$). Now, predict (impute) the averages at the two children regions of our region $R_{l,k}$ by setting

$$\tilde{\beta}_{l+1,2k-1} = Ave\{\pi_{l,k}|R_{l+1,2k-1}\} \qquad \tilde{\beta}_{l+1,2k} = Ave\{\pi_{l,k}|R_{l+1,2k}\}.$$

Note that while the function $\pi_{l,k}$ is global, but it is only fit to the local information at and around the given region, and used to predict the averages only at the children of the given region; a different function $\pi_{l,k}$ is needed for each region $R_{l,k}$. By repeating this procedure for all of the regions at level $l$, we will obtain predicted averages for all of the regions at level $l + 1$.

This two scale refinement scheme allows constructing the average interpolating wavelet transform. Consider the differences between the actual and predicted averages at level $l + 1$ and define the wavelet coefficients

$$\alpha_{l,k} = \frac{Vol(R_{l+1,2k})}{\sqrt{VolAv_{l+1}}} \times (\beta_{l+1,2k} - \tilde{\beta}_{l+1,2k}), \quad l \geq 0, \quad k = 1, ..., 2^l.$$

The normalization in this formula achieves two goals. First, note that there is no need to keep the differences for the other child, $\alpha'_{l,k} = \frac{Vol(R_{l+1,2k-1})}{\sqrt{VolAv_{l+1}}} \times (\beta_{l+1,2k-1} - \tilde{\beta}_{l+1,2k-1})$, due to the redundancy in our pyramid of averages. Indeed, Eq. (2.1) implies that

$$(\beta_{l+1,2k-1} - \tilde{\beta}_{l+1,2k-1})Vol(R_{l+1,2k-1}) + (\beta_{l+1,2k} - \tilde{\beta}_{l+1,2k})Vol(R_{l+1,2k}) = 0,$$

or more clearly, due to our scaling,

$$(2.3) \qquad\qquad \alpha_{l,k} + \alpha'_{l,k} = 0.$$

Second, the normalization coefficient has the right order of magnitude in order to make these wavelet coefficients similar to the usual wavelet coefficients in that the magnitude of "white noise" wavelet coefficients stabilize to have comparable magnitudes across levels.

It can be easily seen that the original function $f$ can be reconstructed given $\beta_{0,1}$ and all of $\alpha_{l,k}$. Indeed, this backward wavelet transform process proceeds as follows. We first use $\beta_{0,1}$ to compute the predicted averages $\tilde{\beta}_{1,1}$ and $\tilde{\beta}_{1,2}$ using the two scale refinement scheme, and then set $\beta_{1,2} = \tilde{\beta}_{1,2} + \alpha_{0,1} \times \frac{\sqrt{VolAv_1}}{Vol(R_{1,2})}$. Now $\beta_{1,1}$ can be uniquely computed using a similar formula and the redundancy relationship of Eq. (2.3) between the two children of each region. This gives us the averages at level $l = 1$. Repeating this process we can obtain the averages at level $l = 2$, then at level $l = 3$, etc. The function $f$ is given as the limit of these averages.

Now, note that in order to instantiate the above construction on manifolds, we need two ingredients: a method for obtaining the dyadic partition $\{R_{l,k}\}$, and a method for constructing the interpolants $\pi_{l,k}$. The next two subsections address these issues.

2.2. **Partitioning.** In the classical setting of the real line [7], generating the dyadic subdivision is trivial. While this is not as easy in the setting of manifolds, we believe that an appropriate procedure can be devised. Consider the general problem of partitioning a given region $\Omega \subseteq M$ into two connected regions. One can pick two points $p_1$ and $p_2$ in $\Omega$ and define the first region $\Omega_1$ as the set of points closer to $p_1$ than to $p_2$ as measured by the shortest distance within $\Omega$ (i.e. the length of the shortest path that remains within $\Omega$); the second region $\Omega_2$ is defined similarly. This initial subdivision can be improved by optimizing the choice of points $p_1$ and $p_2$ in order to maximize some measure of quality, such as the relative sizes of the generated regions.

Note that since the regions in the partitions get smaller and smaller with the increasing level $l$, after some level it will be possible to identify (e.g. via the exponential map) a region with a subset of the tangent space of a point within the region. Thereafter, it will be possible to apply some of the existing Euclidean approaches, such as 2-means (i.e. $k$-means clustering algorithm with $k = 2$), to partition this subset of the tangent space and transfer the partition back to the manifold. We anticipate that some of the results in [16], such as the existence of upper limits on the eccentricities of obtained regions can be transferred into the manifold setting.

Since the main goal of this paper is to obtain a wavelet construction on discrete data, we skip a deeper discussion of this issue in the continuous case, and make an assumption that some "good" subdivision of the manifold is given. In practice, as discussed in the next section, we use a spectral embedding (somewhat similar to the diffusion map) of the manifold into a high-dimensional Euclidean space, and recursively subdivide it using the 2-means algorithm.

### 2.3. Interpolation.

Our interpolant construction is based on the average interpolating variational splines that appear in Pesenson [13]. However, the splines of [13] cannot be used directly here because they do not produce a constant interpolant when all of the prescribed averages are the same. We modify his formulation in order to attain this requirement leading to a vanishing moment, and we show how more vanishing moments can be obtained. We also derive formulas for the interpolant that allow for generalizations and efficient computation in the discrete setting.

Consider the general problem: given disjoint (except for possibly common boundaries) regions $\{\Omega_r\}_{r=1}^m$ on a closed Riemannian manifold $M$, construct a function $\pi : M \to \mathbb{R}$ which attains the prescribed averages,

$$(2.4) \qquad Ave\{\pi|\Omega_r\} = \beta_r, \quad r = 1, ..., m.$$

Let $\Delta$ be the Laplace-Beltrami operator on $M$; define the interpolant $\pi$ as the minimizer (in a suitable space of functions) of the following energy

$$\int_M (\Delta^{p/2}\pi)^2$$

subject to the constraints of Eq. (2.4); here $p$ is an appropriate positive even integer. The objective function is strictly convex over the feasible set, and, so, the problem has a unique solution. Note that if all of $\beta_i$ values are the same, then the constant function $\pi = const$ would be both feasible and make the objective equal to zero, which means that the constant function would be the minimizer, as required.

The choice of the objective function is natural due to the following reasons. First, on the real line, by setting $p = 2$ we obtain the functional $\int(\pi'')^2$ which is minimized by cubic splines. Second, this is a generalization of the Laplacian energy $\int_M (\Delta\pi)^2$ which is well-known to provide a measure of smoothness for functions on manifolds, see e.g. [2]. As a result, one can say that in some sense our optimization problem is finding the smoothest function on the manifold satisfying the averaging constraints.

Next, we show that the solution of this minimization problem can be written in terms of the poly-harmonic Green's function, which leads to an efficient solution in practice. Our discussion is rather informal, as we will only use the resulting formula to guide our construction in the discrete case.

The following preliminaries will be needed. Let $\{\lambda_n, \phi_n\}_{n=0}^\infty$ be the Laplace-Beltrami eigenvalues and eigenfunctions; we assume that an orthonormal set of eigenfunctions was chosen. Note that for $\lambda_0 = 0$ the corresponding eigenfunction $\phi_0 = const = 1/\sqrt{Vol(M)}$ due to the unit $L^2(M)$ normalization. Since all eigenfunctions $\phi_n$ for $n > 0$ are orthogonal to $\phi_0$, we obtain that $\int_M \phi_n = 0$. An $L^2(M)$ function $f : M \to \mathbb{R}$ can be expanded as

$$f = a_0\phi_0 + a_1\phi_1 + a_2\phi_2 + a_3\phi_3 + ...$$

where $a_n = \int_M f\phi_n$. The generalized delta function can be defined by the completeness relation $\delta(x, y) = \sum_{n=0}^{\infty} \phi_n(x)\phi_n(y)$. Note that, as usual, it satisfies

$$\int_M \delta(x, y)f(x)dx = \sum_{n=0}^{\infty} \phi_n(y) \int_M \phi_n(x)f(x)dx = \sum_{n=0}^{\infty} \phi_n(y)a_n = f(y).$$

Consider the $p$-harmonic kernel

$$G(x, y) = \sum_{n=1}^{\infty} \frac{\phi_n(x)\phi_n(y)}{\lambda_n^p}.$$

We will assume that $p$ is chosen such that the series $\sum_{n=1}^{\infty} 1/\lambda_n^p$ converges. If $\Delta_x$ is the Laplace-Beltrami operator acting on the variable $x$, it is easy to verify that

$$
\begin{aligned}
(2.5) \qquad \Delta_x^p G(x, y) & = \sum_{n=1}^{\infty} \Delta_x^p \phi_n(x)\phi_n(y)/\lambda_n^p = \sum_{n=1}^{\infty} \lambda_n^p \phi_n(x)\phi_n(y)/\lambda_n^p \\
& = \sum_{n=1}^{\infty} \phi_n(x)\phi_n(y) = \delta(x, y) - \frac{1}{Vol(M)}.
\end{aligned}
$$

The subtracted term is $\phi_0(x)\phi_0(y) = 1/\sqrt{Vol(M)} \times 1/\sqrt{Vol(M)}$ and is needed as the summation in the delta function starts at $n = 0$. Finally, it is easy to see that $\int_M G(x, y)dx = 0$.

Now we will derive the optimality conditions for our optimization problem. The prescribed average requirements can be rewritten as

$$\int_{\Omega_r} \pi = \beta_r Vol(\Omega_r),$$

which then are incorporated into the problem to get the functional

$$\mathcal{E}(\pi) = \int_M (\Delta^{p/2}\pi)^2 - \sum_{r=1}^{m} c_r(\int_{\Omega_r} \pi - \beta_r Vol(\Omega_r)),$$

where $c_r$ are the Lagrange multipliers (in fact, $-c_r$ are the Lagrange multipliers). The first variation of the functional $\mathcal{E}$ after dropping the second order terms in $\delta\pi$ is

$$\delta\mathcal{E} = \mathcal{E}(\pi + \delta\pi) - \mathcal{E}(\pi) = \int_M 2\Delta^{p/2}\pi\Delta^{p/2}(\delta\pi) - \sum_{r=1}^{m} c_r \int_{\Omega_r} \delta\pi.$$

Now, notice that for any function $f$, it is true that $\int_{\Omega_r} f(x)dx = \int_M f(x) \int_{\Omega_r} \delta(x, y)dydx$, and so we get

$$
\begin{aligned}
\delta\mathcal{E} & = \int_M 2\delta\pi\Delta^p\pi - \int_M \delta\pi \sum_{r=1}^{m} c_r \int_{\Omega_r} \delta(x, y)dydx \\
& = \int_M \delta\pi \left( 2\Delta^p\pi - \sum_{r=1}^{m} c_r \int_{\Omega_r} \delta(x, y)dy \right)
\end{aligned}
$$

Here we also used integration by parts for the first term and the fact that the manifold has no boundary. Setting this equal to zero and remembering that $\delta\pi$ is any small variation, we get that

$$2\Delta^p\pi - \sum_{r=1}^m c_r \int_{\Omega_r} \delta(x,y)dy = 0$$

is satisfied by the optimal function $\pi$.

To solve this, we make the substitution $\pi(x) = \rho(x) + \frac{1}{2}\sum_{r=1}^n c_r \int_{\Omega_r} G(x,y)dy$, where $G(x,y)$ is the $p$-harmonic kernel. Notice that,

$$\Delta_x^p \int_{\Omega_r} G(x,y)dy = \int_{\Omega_r} \Delta_x^p G(x,y)dy = \int_{\Omega_r} \left(\delta(x,y) - \frac{1}{Vol(M)}\right),$$

where we used Eq. (2.5). As a result, we find that $\rho$ satisfies

$$2\Delta^p\rho - \sum_{r=1}^n c_r \frac{Vol(\Omega_r)}{Vol(M)} = 0.$$

Since we assumed that $M$ is without a boundary, this equality can only be satisfied if the constant term in the equation is zero, i.e.

$$(2.6) \qquad \sum_{r=1}^m c_r Vol(\Omega_r) = 0,$$

and then the solution is $\rho = const \triangleq c_0$.

Based on this, we can write the solution of the optimization problem as

$$(2.7) \qquad \pi = c_0 + \sum_{r=1}^m c_r \int_{\Omega_r} G(x,y)dy,$$

where $c_r, r = 0, 1, ..., m$ are some coefficients and $G(x,y)$ is the $p$-harmonic kernel.

There are $m+1$ coefficients in the formula for $\pi$, and to determine them uniquely we will need as many equations. These equations are provided by the $m$ averaging constraints of Eq. (2.4) and the equality of Eq. (2.6). Note that all of these equations are linear. Indeed, let $A$ be the symmetric $m \times m$ matrix with entries

$$A_{rr'} = \int_{\Omega_r} \int_{\Omega_{r'}} G(x,y)dydx,$$

and let $a$ be an $m \times 1$ vector with entries $a_r = Vol(\Omega_r), r = 1, ..., m$. We also define the $m \times 1$ vector $b$ with entries $b_r = \beta_r Vol(\Omega_r), r = 1, ..., m$. The linear system determining the coefficients in the Eq. (2.7) is given by

$$\left[\begin{array}{cc} 0 & a^t \\ a & A \end{array}\right] \left[\begin{array}{c} c_0 \\ c \end{array}\right] = \left[\begin{array}{c} 0 \\ b \end{array}\right],$$

where $c$ is the $m \times 1$ vector with entries $c_r, r = 1, ..., m$.

**Higher vanishing moments.** More vanishing moments can be obtained within this framework as well. Following [5], vanishing moments in the manifold setting are defined as the wavelets being orthogonal to a Laplace-Beltrami eigenfunction. For example, assume that we would like to obtain an average interpolating construction where the wavelets are orthogonal to the eigenfunction $\phi_1$ corresponding to the eigenvalue $\lambda_1$. This means that if the interpolation problem is given averages that arise from the function $\phi_1$, i.e. that if the $\beta_r$ in Eq. (2.4) are given by $\beta_r = Ave\{\phi_1|\Omega_r\}$, then the interpolant must be given by $\pi = \phi_1$.

This can be achieved by defining the interpolant as the minimizer of the following functional:

$$\int_M (\Delta^{p/4}(\Delta - \lambda_1)^{p/4}\pi)^2.$$

Now, functions $\pi \in span\{\phi_0, \phi_1\}$ make this functional equal to zero, and so, when the imposed averages come from a function in $span\{\phi_0, \phi_1\}$, then that function will be uniquely reconstructed. This holds assuming that the number of regions in the problem satisfies $m \geq 2$, and that there is a gap between $\lambda_1$ and the next eigenvalue $\lambda_2$. If the latter is not true, namely if $\lambda_1$ is not a simple eigenvalue, then the number of the additional vanishing moments will be equal to the multiplicity of $\lambda_1$. Accordingly, the number of regions in the problem should not be less than the total number of vanishing moments. [1]

As before, an explicit form for the solution of this optimization problem can be obtained. Without going into the details, we just mention that the solution has the form

$$\pi = c_{-1}\phi_1 + c_0 + \sum_{r=1}^{m} c_r \int_{\Omega_r} G(x, y)dy,$$

where the kernel $G(x, y)$ is defined by

$$G(x, y) = \sum_{n=2}^{\infty} \frac{\phi_n(x)\phi_n(y)}{\lambda_n^{p/2}(\lambda_n - \lambda_1)^{p/2}}.$$

There are $m+2$ coefficients in the formula for $\pi$, and to determine them uniquely we will need as many equations. As before, $m+1$ of these equations are provided by the $m$ averaging constraints of Eq. (2.4) and the equality of Eq. (2.6). The additional equation that must be satisfied is

$$\sum_{r=1}^{m} c_r \int_{\Omega_r} \phi_1 = 0.$$

A straightforward generalization of this technique can be used to obtain further vanishing moments.

## 3. Discrete Setting

3.1. **Discrete Setup.** In order to formulate our construction in the most general setting, we consider a simple connected weighted graph $G$ with vertex set $V$ of size $N$; in the following, we identify $V$ with $\{1, 2, ..., N\}$. We assume that non-negative weights are associated with edges, and, perhaps, with vertices. For a real valued function on the graph vertices, $f : V \to \mathbb{R}$, we will denote by $f(i)$ its value on the $i$-th vertex. We sometimes will think of such functions as vectors in $\mathbb{R}^N$.

Using the edge weights, we can construct the symmetric $N \times N$ matrix $W$ which has entry $W_{ij}$ equal to the weight associated with the edge $(i, j)$ connecting the vertices $i$ and $j$; we set $W_{ij} = 0$ when $i$ and $j$ are not connected with an edge (since there are no self-loops, we always have $W_{ii} = 0$). We also construct the diagonal $N \times N$ matrix $S$ which captures the weights associated with the vertices: we set $S_{ii}$ equal to the weight of the $i$-th vertex in the graph. If the graph has no vertex weights, we can simply set $S_{ii} = 1$ for all $i = 1, ..., N$; note that there are other options, one of which is setting $S_{ii} = \sum_j W_{ij}$, the sum of all edge weights incident to the $i$-th vertex.

---

[1] In the average interpolation construction the required number of neighboring regions can be achieved by considering higher order neighborhoods. In addition, in the initial level of the dyadic partitioning of the manifold ($l = 0$) one needs to start with at least as many partitions as the number of the vanishing moments.

The role that the Laplace-Beltrami operator plays in the continuous case will be played in the discrete setting by the graph Laplacian. Let $D$ be the diagonal matrix with entries $D_{ii} = \sum_j W_{ij}$. The graph Laplacian $L$ is defined as

$$L = S^{-1}(D - W).$$

While our problem setup started with a weighted graph and arrived to the Laplacian matrix $L$, our construction can also be applied when one starts with the Laplacian matrix $L$ and infers from it the weighted graph. This is a natural way of dealing with point clouds sampled from a low-dimensional manifold, a setting common in manifold learning. There are numerous ways for computing Laplacians on point clouds, see [1, 4, 3, 6]; most of them fit into the above form $L = S^{-1}(D - W)$, and so, they can be used to infer a weighted graph. Some of the point cloud Laplacians converge to the Laplace-Beltrami operator of the underlying manifold as the density of point sampling is increased, a fact that makes the parallels between our discrete and and continuous constructions more apparent.

3.2. **Preliminaries.** Given the Laplacian $L = S^{-1}(D - W)$, we will need to compute its eigenvalues and eigenvectors. Note that $L$ may not be symmetric due to the factor of $S^{-1}$. To avoid solving a non-symmetric eigenvalue problem, we solve the following symmetric generalized eigenvalue problem

$$(D - W)\phi = \lambda S\phi.$$

Note that the eigenvalues are non-negative, with $0 = \lambda_0 < \lambda_1 \leq \lambda_2 \leq \lambda_3 \leq ... \leq \lambda_{N-1}$. As in the continuous case, we can pick an orthonormal set of eigenvectors $\phi_n$, $n = 0, 1, ..., N - 1$, with $\phi_0 = const$. Importantly, these eigenvectors will be orthonormal with respect to the $S$-inner product. Namely, we will have

$$(3.1) \qquad \phi_n S \phi_{n'} = \delta_{nn'},$$

where $\delta_{nn'}$ is Kronecker's delta.

An interesting ramification of Eq. (3.1) relates to how one converts integrals appearing in the continuous setting to the current discrete setting. Indeed, comparing Eq. (3.1) to the continuous orthonormality relation, $\int_M \phi_n \phi_{n'} = \delta_{nn'}$, it becomes apparent that integrals should become weighted sums,

$$(3.2) \qquad \int_M f \to \sum_i S_{ii} f(i).$$

Accordingly, for a function $f : V \to R$, we define its average value over a subset of vertices $R \subseteq V$ as the weighted average:

$$Ave\{f|R\} = \frac{\sum_{i \in R} S_{ii} f(i)}{\sum_{i \in R} S_{ii}}.$$

Based on this formula, it is natural to define the volume of $R$ by $Vol(R) = \sum_{i \in R} S_{ii}$.

### 3.3. Embedding.

Our construction will interact with the graph structure through its embedding into a high dimensional Euclidean space. Consider the map $\mathcal{H} : V \to R^{N-1}$ defined by the formula

$$(3.3) \qquad \mathcal{H}(i) = \left( \frac{\phi_1(i)}{\lambda_1^{p/2}}, \frac{\phi_2(i)}{\lambda_2^{p/2}}, ..., \frac{\phi_{N-1}(i)}{\lambda_{N-1}^{p/2}} \right),$$

to which we will refer as the $p$-harmonic embedding of the graph. Note that this is similar to the Diffusion Map [4] with the difference being in how the individual eigenvectors are scaled, by a factor of $\lambda_n^{-p/2}$ here, versus $e^{-t\lambda_n}$ for the Diffusion Map. When $p = 2$, we obtain the biharmonic embedding introduced in the case of the surface meshes in [11].

The relevance of the $p$-harmonic embedding becomes clear due to the following connection. We will need a discrete equivalent of the $p$-harmonic Green's kernel, which one naturally defines similarly to the continuous case as:

$$G(i,j) = \sum_{n=1}^{N-1} \frac{\phi_n(i)\phi_n(j)}{\lambda_n^p}.$$

It is easy to see that the following holds,

$$G(i,j) = \mathcal{H}(i) \cdot \mathcal{H}(j),$$

where the dot signifies the usual dot product of vectors.

In practice two issues need to be addressed. First, for efficiency reasons, it is impractical to compute the complete eigen-decomposition of the Laplacian, and so we truncate the embedding to only include the eigenvectors corresponding to the smallest few-hundred eigenvalues. Second, a choice of $p$ should be made. While one can use this as a tuning parameter, we recommend setting $p$ in such a way that the denominator in the embedding of Eq. (3.3) grows at least linearly in $n$. In the point cloud case where the discrete Laplacian approximates the Laplace-Beltrami operator of a smooth manifold $M$, we know by the Weyl estimate that $\lambda_n \sim n^{2/dim(M)}$. Thus, to make the denominator of at least linear order one needs to set $p \geq dim(M)$. Various techniques for estimating the dimensionality of a point cloud are available, for our purposes a crude estimate (perhaps, based on fitting a line to the log-log scatter-plot of $(n, \lambda_n)$ on the plane) will suffice. Setting the value of $p$ too high, will effectively discard the higher eigenvectors, and will result in numerical stability issues at the very fine levels, as there will not be enough variation in the kernel to provide enough degrees of freedom for interpolation.

Finally, we note that a straightforward generalization of our construction comes by using other embeddings instead of the embedding of Eq. (3.3). Our specific choice was driven by the formula for the interpolant in the continuous case which requires the $p$-harmonic Green's kernel. This kernel, in turn, arises due to our choice of the functional to minimize. If a different embedding is used, implicitly the minimized functional will assume another form. However, in practice, the explicit form of the functional is not needed in order to complete the construction, the only entity that our construction interacts with is the embedding one chooses to use.

### 3.4. Partitioning.

As in the continuous case, in order to realize the average-interpolating wavelet scheme we first need to construct a dyadic partitioning of the graph $G$. While a variety of graph partitioning techniques exist, we base our partitioning on the spectral clustering algorithm of [12] as applied to the embedding of the Eq. (3.3). To obtain a binary tree of partitions, we start with the graph itself as the root. At every step, a given region (a subset of the vertex set) of graph $G$ is split into two children partitions by running the 2-means clustering algorithm ($k$-means with $k = 2$) on the embedding of Eq. (3.3) restricted to the vertices of the given partition. This process is continued in recursion at every obtained region.

In practice two issues need to be addressed. First, due to the finite size of the graph, we cannot continue partitioning the graph ad infinitum. The maximum number of levels, $l_{max}$, should be chosen based on the number of graph vertices. We usually set the number of levels as $l_{max} = \lfloor \log_2 N \rfloor$; smaller values can be set if efficiency is a consideration or if a less detailed analysis of data is needed. Second, we did not let the non-leaf regions (i.e. regions at levels $l < l_{max}$) to become less than a certain size (2 vertices in our experiments). As a result, we would get nodes in the binary tree which have a single child (i.e. a region would be small enough so that it is not split anymore). On the other hand, the regions at level $l_{max}$ would each contain a single node of the graph. Thus, in reality, our tree is not strictly a binary tree at the level $l_{max}$. The redundancy in the average-interpolating pyramid for a binary partition tree allows one to keep one wavelet coefficient at each parent node. To accommodate non-binary trees in our implementation, we keep the wavelet coefficients at the children nodes (namely, we keep $\alpha_{l,k}$ at node $R_{l+1,2k}$, then $\alpha'_{l,k}$ at $R_{l+1,2k-1}$, and so on if there are more children) which makes our representation redundant. When modifying the wavelet coefficients (e.g. in order to achieve denoising), one has to make sure that the straightforward generalization to non-binary trees of the redundancy relationship of Eq. (2.3) holds.

3.5. **Interpolation.** The second ingredient in the average-interpolating wavelet scheme is interpolation. To achieve interpolation on graphs, we simply modify the formulas provided in the continuous case. First, the $p$-harmonic Green's kernel $G(i, j)$ is computed using the truncated $p$-harmonic embedding described above. Now, the interpolant is given by the formula of Eq. (2.7) where the coefficients are computed using the same linear system of equations as in the continuous case. When setting up this system of equations, we replace all of the integrals by weighted sums according to the rule in Eq. (3.2).

The use of the truncated embedding can be interpreted as follows. Suppose that in the embedding we are using only the eigenvectors up to index $n_{max}$. One can prove that the obtained interpolant gives the minimizer of the functional $\int_M (\Delta^{p/2} \pi)^2$ among $\pi \in span\{\phi_0, \phi_1, ..., \phi_{n_{max}}\}$, subject to the averaging constraints. This means that our interpolant $\pi$ is a linear combination

$$\pi = a_0 + a_1 \phi_1 + a_2 \phi_2 + ... + a_{n_{max}} \phi_{n_{max}},$$

and among all such linear combinations that satisfy the averaging constraints, we choose the one that minimizes

$$(3.4) \qquad \int_M (\Delta^{p/2} \pi)^2 = \sum_{n=1}^{n_{max}} a_n \lambda_n^p.$$

While it seems that one could directly solve this optimization problem without any recourse to the solution formula of Eq. (2.7), in practice we found that the solutions obtained using such direct methods may become numerically unstable and lead to discontinuities in the scaling functions [2]. On the other hand, this formulation of the problem bears more resemblance to the classical approach of Donoho [7]. Indeed, in the setting of the real line, the solution is sought as a linear combination of the monomials, and the averaging constraints are satisfied by a unique linear combination of the monomials. In our setting, monomials are replaced by Laplace-Beltrami eigenfunctions, and due to the variable number of neighbors there will be many linear combinations that satisfy the averaging constraints. To single out a linear combination among them we pick one that minimizes a measure of smoothness as given by the expression in Eq. (3.4).

One issue that arises during interpolation when dealing with point clouds is the question of which regions should be considered neighbors. The simplest approach would be to declare two regions neighbors if there is an edge connecting a vertex from one region to a vertex from the other. Note that the graph structure for point clouds is inferred from the point cloud's Laplacian. These Laplacians may have rather large stencils: for example, in one approach, the stencil at a given point contains a fixed number, sometimes 50-100, of its nearest neighbor points. If this is the case, then even a region containing two points will have a large number of neighboring regions, which can be counter-intuitive to the notion that the point cloud is a representation of an underlying low-dimensional manifold; this will also prevent the wavelets from having supports smaller than a certain size. To circumvent this issue, we look at the total weight sum for all edges running between the two regions, and if this "cut size" is large enough, then we declare the two regions to be neighbors. For each region, we set the minimum threshold on the cut size equal to a fraction of the largest cut size at the given region.

## 4. Numerical Examples

In order to investigate the empirical properties of the average interpolating construction introduced in this paper, we ran a set of experiments using a number of synthetic and real data sets. For the experiments, we compute and visualize the scaling functions at different levels $l$. There is one scaling function for each region $R_{l,k}$, which can be obtained by setting the function average at that region $R_{l,k}$ equal to one, and then applying the refinement process. To get scaling functions supported at similar locations, we pick a point on the graph, and at each level show the scaling function corresponding to the region that contains the picked point. For the ease of visualization, all functions are scaled to have the maximum absolute value of 1 and the following color coding is used: dark red represents larger positive values, neutral green corresponds to zero, and dark blue represents larger (in absolute value) negative values. Wavelets are not shown because the wavelet at a given parent region is simply the weighted difference between the scaling functions corresponding to the two child regions. In all of the experiments we used the value of $p = 2$. In all of the examples, except the one-dimensional case, we only look at the immediate neighbors of a region, out of which we retain only the ones that have large enough cut size (at least $1/8$ of the largest cut size at the given region). All presented example are based on a single vanishing moment.

---

[2]For example, we found the following direct approach to lead to discontinuities. First, eliminate $a_0$ using one of the averaging constraints. Next, rescale the unknowns so that the expression that is minimized is simply the sum of squares of the unknowns ($b_n \triangleq a_n \lambda_n^{p/2}$), and then use the pseudo-inverse to solve the remaining averaging constraint equations.

Figure 1 shows the scaling functions in one-dimensional case. The discrete Laplacian $L$ is obtained using the three point stencil $[-1, 2, -1]$. Our subdivision process reproduces the dyadic subdivision of the interval up to a small numerical error. Different scaling functions result based on how many neighbors of a given interval are used during interpolation. For example, if one uses the immediate neighbors of an interval only, the function shown in Figure 1 (a) results; but if both the immediate neighbors and their neighbors are used, then Figure 1 (b) results. In the original construction of Donoho, the first case would require a quadratic polynomial, and the second case – a quartic polynomial, and so on. In our construction, independently of the number of neighbors, the form of the interpolant stays the same. As a result, except for Figure 1 (a), our scaling functions are different from Donoho's, which can be seen by examining the empirical two scale relation coefficients of our construction. Note that as in the original construction, increasing the number of neighbors results in smoother but less localized scaling functions.

Figure 2 (a) depicts a graph representing the road network for Minnesota, with edges showing the major roads and vertices being their intersections. Every edge has unit weight, while every vertex's weight is set equal to its degree; as a result, the Laplacian $L$ equals to the normalized graph Laplacian. The obtained scaling functions at different levels are shown in Figures 2 (b-f).

Figure 3 shows the scaling functions on the "Swiss Roll", a 2-dimensional manifold with boundary. Here, it is represented as a uniformly sampled point cloud and the Laplacian is computed using the point cloud Laplacian of [4]. Since the ingredients used in our construction are intrinsic, the scaling functions are supported on contiguous regions along the manifold.

Figure 4 depicts the scaling functions on a planar domain with irregular boundary. Once again, this is considered as a point cloud and the Laplacian of [4] is used. Note that the supports of scaling functions adapt to the shape of the domain by initially spreading within the "bulb" before starting to grow into the narrow region connecting the "bulbs". This behavior is similar to how spectral distances, such as the diffusion distance, behave.

In Figure 5 we exemplify a straightforward application of our construction to smoothing a noisy function on a manifold. We added Gaussian noise to the smooth function on the sphere in Figure 5 (a), to obtain the function in Figure 5 (b) with the SNR of 1.24 dB. We compute the wavelet transform of the function in (b) and then set all of the wavelet coefficients after certain level ($l \geq 5$) equal to zero. Inverting the transform, we obtain the function depicted in 5 (c) with the SNR of 18.2 dB. Clearly, more sophisticated wavelet denoising approaches can be applied in general.

Another simple application of our construction is to the problem of regression on manifolds. For a function given on the S-manifold in Figure 6 (a), we randomly sample 75 points shown in Figure 6 (b). The function values at these sampled points are used to reconstruct the original function on the entire manifold using the following procedure, similar in the spirit to that of [8]. First, whenever a region contains at least one of the sampled points where the function value is provided, we can estimate the function's average at that region using these sample values (in the obvious way). These averages are put in the pyramid, and then we compute the wavelet coefficients $\alpha_{l,k} \sim (\beta_{l+1,2k} - \tilde{\beta}_{l+1,2k})$ and $\alpha'_{l,k} \sim (\beta_{l+1,2k-1} - \tilde{\beta}_{l+1,2k-1})$ whenever all of the involved terms can be estimated. This requires that the children regions $R_{l+1,2k}$ and $R_{l+1,2k-1}$, the parent region $R_{l,k}$, and the neighbors of the parent used for interpolation all contain at least a single sampled point with a known function value. All of the wavelet coefficients that cannot be estimated are set to zero. In the final step, we need to reconcile the values of $\alpha_{l,k}$ and $\alpha'_{l,k}$ so that the redundancy relationship $\alpha_{l,k} + \alpha'_{l,k} = 0$ is satisfied in the entire pyramid; this is achieved by setting $\alpha_{l,k} = (\alpha_{l,k} - \alpha'_{l,k})/2$ and $\alpha'_{l,k} = (\alpha'_{l,k} - \alpha_{l,k})/2$. The reconciliation step is need because the sampled points may distribute unevenly between the regions, and so the *estimated* averages in the pyramid may not satisfy the redundancy relationship of Eq. (2.1) that would have been satisfied for the *actual* averages. Next, running the backward wavelet transform, we obtain the reconstructed function. In the example provided here, this process yields the function shown in Figure 6 (c); the normalized RMSE of this reconstruction is 2.83%.

## 5. Summary and Future Work

We introduced a multiscale construction based on a generalization of the average interpolating scheme of Donoho [7]. The construction yields wavelets suitable for analyzing functions on graphs and point clouds. Our main goal in this paper was to call the reader's attention to the possibility of such a construction by motivating it in the continuous case and adapting to the discrete case. As a result, this work provides only a small, first step and has limitations that suggest topics for future work. A first topic suitable for further investigation is to more formally characterize the theoretical properties of the construction in the continuous case. For example, one observes experimentally that our refinement process converges to yield continuous scaling functions. Providing a rigorous proof of this observation would be interesting. A second topic for future research would be to develop theoretical bounds on the decay of the wavelet coefficients in the discrete case. In fact, it seems that some of the results from [8] should easily generalize to our setting. Finally, it would be interesting to investigate applications of the construction to areas such as machine learning, image processing, and computer graphics.

## Acknowledgments

This project would not have been possible without the help of people who provided the software and data. We are grateful to Mauro Maggioni for making publicly available the "Diffusion Geometry", "Diffusion Wavelets", and "Spectral Graph Theory Demo" libraries. We thank David Hammond for distributing the "SGWT: Spectral Graph Wavelets Toolbox". The Minnesota traffic graph is from `matlab_bgl` software package by David Gleich. We thank Michael Chen for making publicly available the `litekmeans` program. The wonderful course notes [15] were extremely helpful as an introduction to the subject.

Department of Mathematics and Computer Science, Drew University, 36 Madison Ave, Madison NJ 08854 USA

*E-mail address*: `rrustamov@drew.edu`

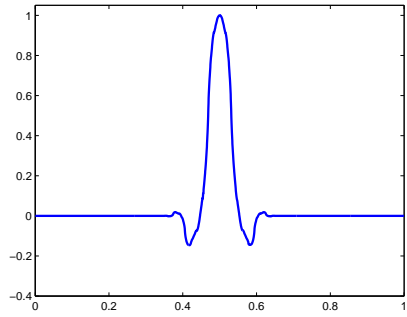

(a) First order neighborhoods

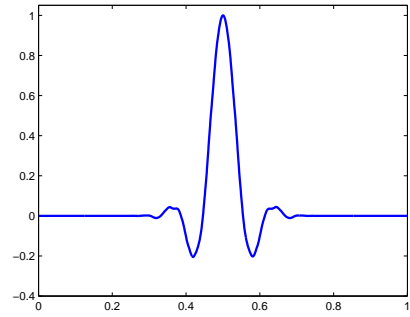

(b) Second order neighborhoods

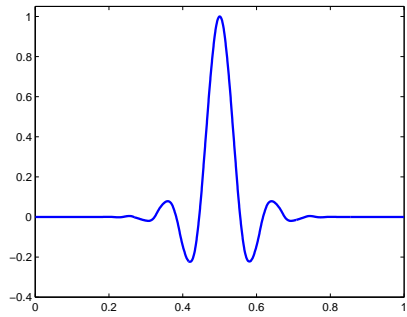

(c) Third order neighborhoods

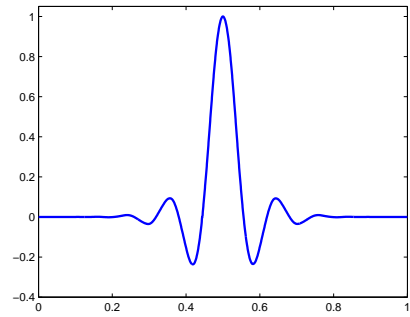

(d) Fourth order neighborhoods

FIGURE 1. Scaling functions on the real line at a single level (scale) for different choices of the neighborhood order.

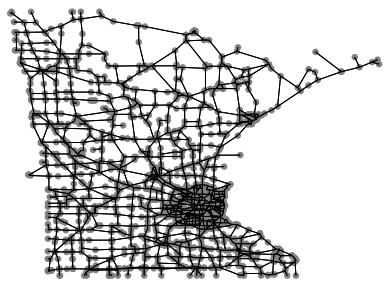

(a) Minnesota road graph

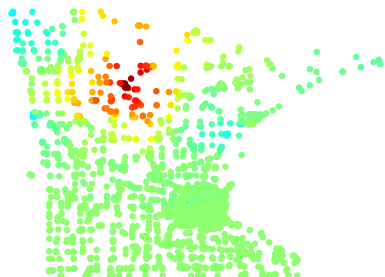

(b) A scaling function at level $l = 4$

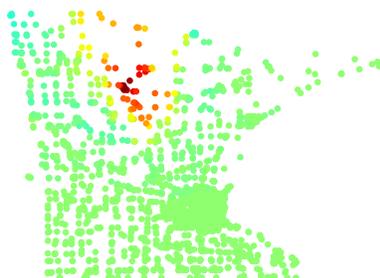

(c) A scaling function at level $l = 5$

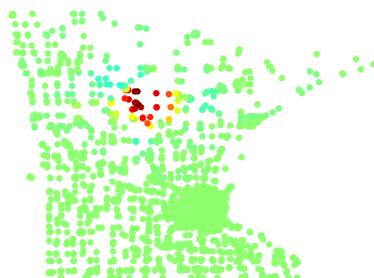

(d) A scaling function at level $l = 6$

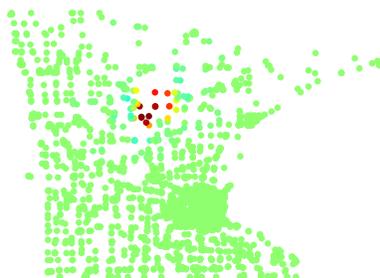

(e) A scaling function at level $l = 7$

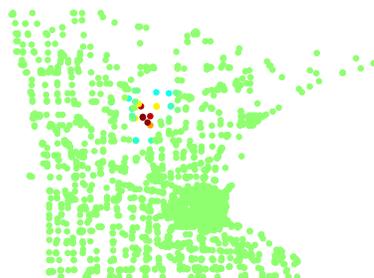

(f) A scaling function at level $l = 8$

FIGURE 2. Scaling functions at different levels on the Minnesota road graph.

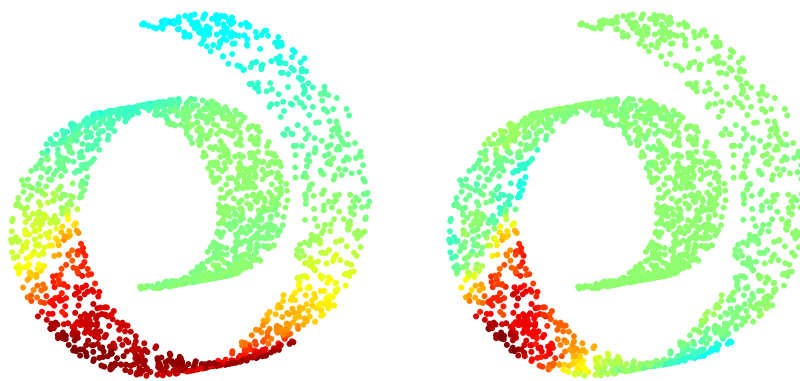

(a) A scaling function at level $l = 2$     (b) A scaling function at level $l = 3$

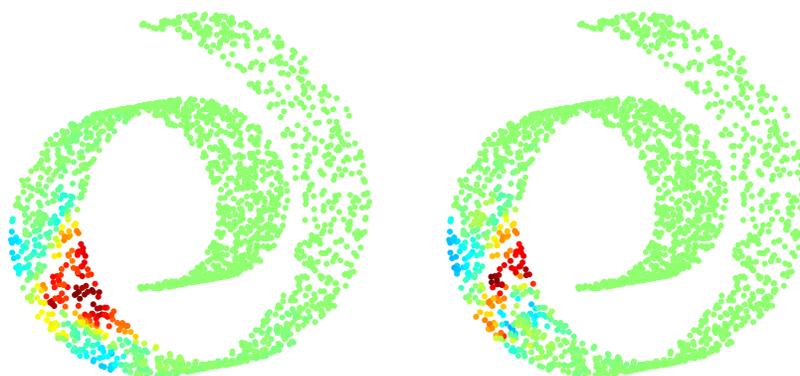

(c) A scaling function at level $l = 4$     (d) A scaling function at level $l = 5$

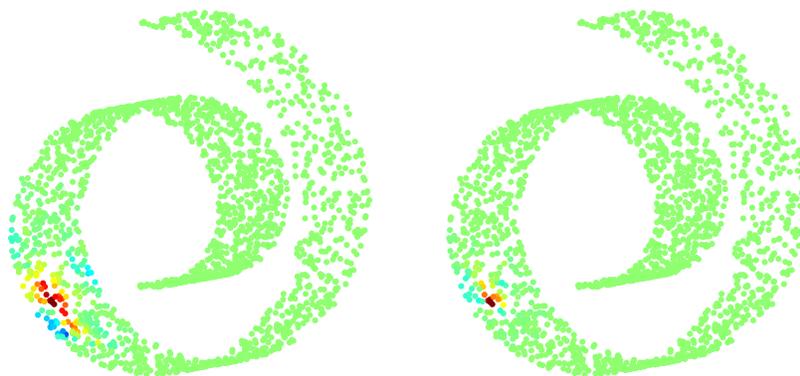

(e) A scaling function at level $l = 6$     (f) A scaling function at level $l = 7$

FIGURE 3. Scaling functions at different levels on the Swiss Roll.

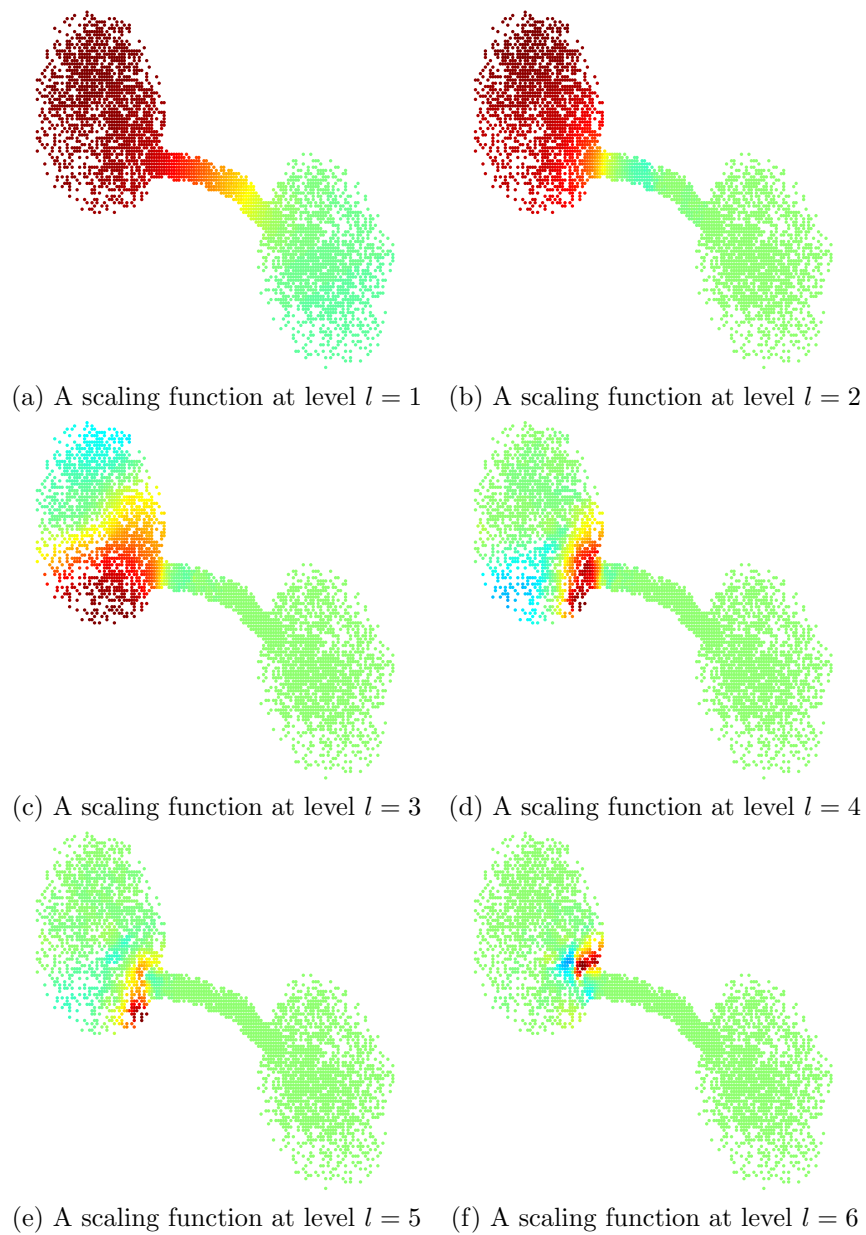

(a) A scaling function at level $l = 1$  (b) A scaling function at level $l = 2$

(c) A scaling function at level $l = 3$  (d) A scaling function at level $l = 4$

(e) A scaling function at level $l = 5$  (f) A scaling function at level $l = 6$

FIGURE 4. Scaling functions at different levels on an irregular planar domain.

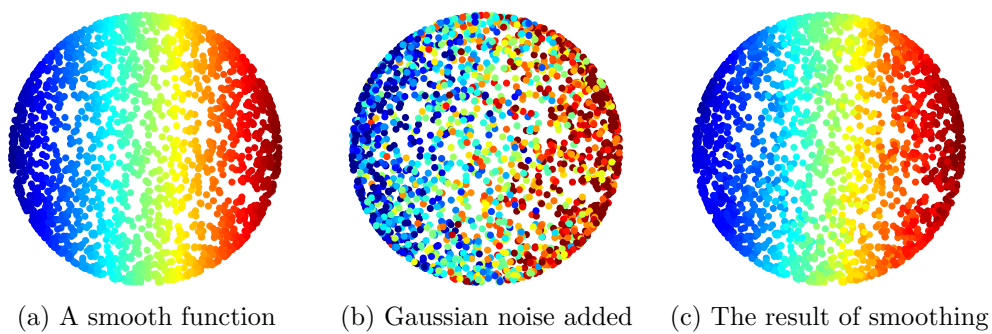

(a) A smooth function  (b) Gaussian noise added  (c) The result of smoothing

FIGURE 5. Smoothing out white noise using wavelets on the sphere.

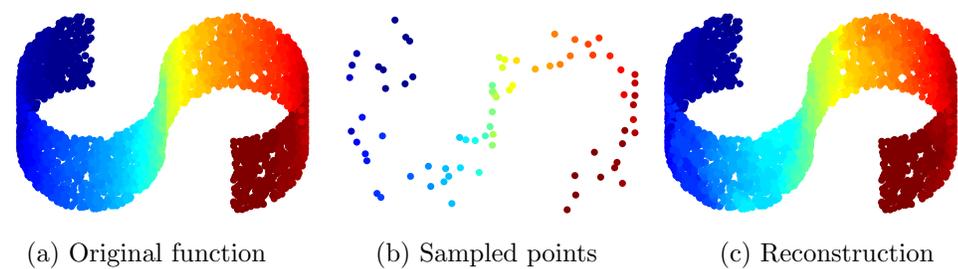

(a) Original function  (b) Sampled points  (c) Reconstruction

FIGURE 6. Reconstructing a function from its values at a subset of points on the S-manifold.